\documentclass[12pt]{article}
\usepackage{makeidx}
\usepackage{amssymb}
 \usepackage{latexsym}
\usepackage{tikz}
\usepackage{pgfplots}
\usepackage{pgfkeys}%
\usetikzlibrary{plotmarks}
\usetikzlibrary{calc}
 \usepackage{array}
\usepackage{arabtex}
 \usepackage[francais]{babel}
 \usepackage[ansinew]{inputenc}
 
 \usepackage{tabularx}
\usepackage{epsfig}
\usepackage{array}
\usepackage[ansinew]{inputenc}
\usepackage{amsfonts,amsmath}
\usepackage{latexsym}
\usepackage{graphics,color,graphicx,shortvrb}

 \usepackage{tabularx}
\newtheorem{Theorem}{Theorem}[section]
\newtheorem{Definition}[Theorem]{Definition}

\newtheorem{Lemma}[Theorem]{Lemma}

\definecolor{ProcessBlue}{cmyk}{1,0,0,0.25}
\definecolor{Black}{cmyk}{0,0,0,1}
\definecolor{Red}{cmyk}{0,1,1,0}
\definecolor{Green}{cmyk}{0.9,0,1,0}
\definecolor{Orange}{cmyk}{0,0.61,0.87,0.1}
\definecolor{Fuchsia}{cmyk}{0.47,0.91,0,0.06}
\definecolor{PineGreen}{cmyk}{0.92,0,0.59,0.25}

\usepackage[plainpages=false,hyperindex=false]{hyperref}
\hypersetup{pdfpagemode=FullScreen}

\begin{document}
 \sloppy
 \begin{center} 
{\Large \bf Transports Regulators of  Networks with Junctions 
Detected by Durations Functions}\\ \mbox{}\\ Jean-Pierre 
Aubin\footnote{VIMADES (Viabilité, Marchés, Automatique, 
Décisions), 14, rue Domat, 75005, Paris, France \\
 aubin.jp@gmail.com, \url{http://vimades.com/aubin/}}\footnote{\textbf{Acknowledgments} 
\emph{This work was partially supported by the Commission of the 
European Communities under the 7th Framework Programme Marie 
Curie Initial Training Network (FP7-PEOPLE-2010-ITN),  project 
SADCO, contract number 264735.}} 
\\ \today

\end{center}
 
\begin{center} 
\textbf{Abstract}
\end{center} 

\emph{This study advocates a mathematical framework   of 
``transport relations'' on a network.  They single out  a subset 
of ``traffic states'' described by  time, duration, position  and 
other traffic attributes (called ``monads'' for short). Duration 
evolutions are non negative, decreasing toward $0$ for incoming 
durations,  increasing from $0$ for outgoing durations, allowing 
the detection  of ``junction states'' defined as  traffic states 
with ``zero duration''. A ``junction relation'' (cross-roads, 
synapses, clearing houses, etc.)  is a 
subset of the transport relation made of junction states. \\
The objective is to construct a ``transport regulator'' 
associating with  traffic states a set of ``celerities'' that 
mobiles circulating in the network can use as velocities. In 
other word,  a network is regarded as a ``provider of  velocity 
information'' to the mobiles for travelling from one departure 
state to an arrival state across a junction relation (a kind of 
geodesic problem). 
\\This investigation assumes that a  system  governs the 
evolution of monads in function of time, duration and position 
using celerities as controls and provide the transport regulator, 
a feedback from transport states to celerities. 
 \\The proposed mathematical framework 
can acclimate road or aerial networks, endocrine (hormonal) or 
synaptic (neurotransmitters) networks, financial or economic 
networks, which motivated this investigation. This framework 
could  probably accommodate computer and even social networks. 
\\This investigation is restricted to junctions between two routes. An 
extension to junctions between a set or prejunction routes and 
postjunction is under way, as well as the development of specific 
examples.}
 
 \mbox{}

 

\section{Introduction}

A \emph{``route''} («~voie~» in French) is a subset $K \subset   
\mathbb{R}^{p} $ of positions $p \in K$, on which travel 
continuously    \emph{``mobiles''} $t \mapsto \xi(t) \in K$. A 
\index{traffic network} \emph{traffic network} is a (finite) 
family of routes intersecting at certain positions, the 
\index{junction} \emph{junctions} of the network, at which some 
traffic attributes of the traffic may be \emph{discontinuous} 
(stopping the mobiles, creating or destroying some attributes, 
for instance). This means that their velocity may be infinite: 
infinite velocities are  called \index{impulse} \emph{impulses}. 
We are investigating the (continuous) circulation of mobiles on 
the routes of the network \emph{punctuated}\footnote{This term is 
borrowed to paleontologists  \glossary{Eldredge (Nils) [1943-]} 
\emph{Nils Eldredge} et \glossary{Gould (Stephen J.) [1941-2002]} 
\emph{Stephen Gould} who proposed it in 1972 for describing 
biological evolution discontinuous at ``punctuated equilibria'' 
(see  the 1977 \cite[\emph{Phylogenesis and 
Ontogenesis}]{Gould77} by \emph{Stephen J. Gould}).} by 
  ``junctions'' between two \emph{prejunction} and \emph{postjunction} 
dates, which, when they are equal, are (instantaneous) 
\emph{impulse junction} dates and positions. 

For instance,  vehicles in road networks evolve along the roads 
through the  cross-roads, where they may be required to stop,  
nervous influx along neurons propagate until they release or 
destroy neurotransmitters in the synapses\footnote{Synapses  are 
``junction'' in Greek, which join  the axon of a 
``pre-synaptic''   neuron and the dendrite of a  
``post-synaptic'', have been discovered in 1897 by 
\glossary{Sherrington (sir Charles Scott) [1857-1952]} 
\emph{Charles Sherrington}.} of a neural network, 
hormones\footnote{from the Greek \emph{horme}, meaning 
\emph{impulse}, \emph{impetus}, have been coined in 1905 by \ 
glossary {Starling (Ernest) [1866-927]}  \emph{Ernest Starling}, 
who isolated them in 1902 together with \glossary {Bayliss 
(William Maddock) [1880-1924]}  \emph{William Bayliss}.} 
circulated in endocrine systems between endocrinal glands and 
receptors,  money   transfers propagate  in  economic-financial 
networks of banks and their exchanges are settled in the clearing 
houses,  information is propagated in a network of computers 
devices until bits\footnote{A bit (contraction of  ``binary 
digit'' proposed by \glossary{Tukey (John Wilder) [1915-2000]} 
\emph{John W. Tukey} in 1947 and next popularized by 
\glossary{Shannon (Claude Elwood) [1916-2001]} \emph{Claude E. 
Shannon} in 1948) is   transmitted by serial or parallel 
transmission one at a time in computing devices.  The encoding of 
data by discrete bits was used in Bacon's cipher (1626), in the 
punched cards  invented by Basile Bouchon and Jean-Baptiste 
Falcon (1732) and developed by Joseph Marie Jacquard (1804), the 
beginning of a long history.} are transferred in their computer 
nodes, etc.

The routes are subsets of positions in a finite dimensional 
vector spaces $\mathbb{R}^{p}$, such as $ \mathbb{R}^{2}$  (or $ 
\mathbb{R}^{3}$ whenever altitude matters) for road networks,   
$\mathbb{R}^{3} $ for aerial trafic and neural networks, red 
globules in blood vessels,  $\mathbb{R}^{p}$ for banks, etc., 
which are associated with adequate \index{traffic attributes} 
\emph{traffic attributes} $x \in \mathbb{R}^{m}$, such as 
celerities (velocities advised to mobiles, measures of  traffic 
jam, number of neurotransmitters, of units of numéraire, etc).

The broad question we ask at a level of abstraction high enough 
to cover these examples is how to regulate the   circulation of 
mobiles through the routes and the junctions of the network.

We impose some requirements  on the traffic attributes at each 
time and each positions of the network which have to be satisfied 
(velocity, jam, etc.).  

The question we would like to answer is:  \emph{how to equip a 
network with a traffic regulator computing and providing   at 
each position and at each time some information transmitted to 
the mobiles travelling along the routes   should or must respect 
to satisfy the requirements?}

\mbox{}

The aim of this study is to suggest a mathematical framework to 
describe a traffic network as an \emph{information provider}. It 
is primarily motivated by road networks\footnote{See \cite[Aubin, 
Bayen \& Saint-Pierre]{absp08hj}, 
\cite[Aubin]{AUB-Leit09,Lax-Hopf}, Chapter~14, p. 563, of 
\emph{Viability Theory.  New Directions}, \cite[Aubin, Bayen \& 
Saint-Pierre]{absp}, devoted to regulation of traffic, and the 
forthcoming \emph{Mathematical Approaches to Traffic Management}, 
\cite[Aubin \& Désilles]{ADA-Traffic}. We also refer to 
\cite[Aubin \& Martin]{Aub-Mart08} for a ``microscopic'' analysis 
of traffic management studying the evolution of mobiles on a 
network without using celerities.}, and, in a less extent by 
neural networks\footnote{See Chapter~8, p.139,  of \mbox{Neural 
Networks and Qualitative Physics: a Viability Approach}, 
\cite[Aubin]{a92ia},  Chapter~7, p. 531, of \emph{La  mort  du  
devin, l'émergence du  démiurge. Essai sur la contingence, la 
viabilité et l'inertie des systèmes}, \cite[Aubin]{mded}, 
\cite[Aubin]{Trend}.} and economic-financial 
networks\footnote{See for instance 
\cite[Aubin]{aconcom98,a01ctfcr} in a  in a connectionnist 
perspective.}. The word ``transport'' being polysemous, 
particularly in mathematics, it is used in many different 
perspectives, such as the Fokker-Plant partial differential 
equations, involving drift and diffusion. This is not this 
meaning which we use in this study. Neither do we address in this 
study the optimal network and transport problems, which have been 
the topic of an extensive literature. For instance,  see 
\emph{Network flows and monotropic 
optimization},\cite[Rockafellar]{RockNetw}, by 
\glossary{Rockafellar (Terry) [1935-]} \emph{Terry Rockafellar} 
and, in the Monge-Kantorovitch perspective, the monographs 
\emph{Optimal Transport: Old and New},\cite[Villani]{Villani} and 
\emph{Optimal Transportation and  
Applications},\cite[Villani]{Villani2} by \glossary{Villani 
(Cédric) [1973-]} \emph{Cédric Villani}. We did not dare to 
introduce ``peregrination'' as a substitute for transport, since 
giving a new meaning to monad was enough to fight polysemy.

Adaptations of this study at a high abstract level providing a 
universal point of view to   potential specific problems  
(traffic engineering, neural network, economic and financial 
networks, etc.) is postponed to future investigations, since they 
require  the contributions of their specialists. Their validation 
to other questions in other fields remains to be done.

\subsection*{Traffic Regulators of the Circulation of Mobiles}

The question we shall try to answer is, in the last analysis, to 
regard a network on which mobiles circulate as \emph{a ``traffic 
regulator'' providing mobiles the velocity that they should or 
must use for guaranteeing the viability and/or the optimality of 
the traffic as well as other requirement}. We call this 
information provided by the traffic network \emph{``celerity''} 
$r(t,p)$ at time $t$ and at position $p$.

More precisely, whenever at time $t$, a given mobile is at 
position $p=\xi(t)$, then its velocity $\xi'(t)$ should coincide 
with the celerity $r(t,p)=r(t,\xi(t))$ provided by the traffic 
regulator. Hence, celerity \emph{``decentralizes''}, so to speak, 
the traffic information used by any mobile\footnote{as prices are 
supposed to do for regulating transactions in theoretical 
economics.}:

\begin{Definition} 
\symbol{91}\textbf{Traffic 
Regulator}\symbol{93}\label{d:TrafficRegulator}\index{} The 
ultimate objective of \emph{traffic regulation}  on a route or a 
network of routes  is to provide at each instant $t$ and at each 
position $p \in K$ on the network a \index{celerity} 
\emph{``celerity''} $r(t,p)$, feeding back the velocity  
$\xi'(t)=r(t,\xi(t))=r(t,p)$ governing the evolution of any 
mobile $t \mapsto \xi(t)$ passing through $\xi(t)=p$ at position 
$p $ at time $t$. Such a map $(t,x) \mapsto r(t,x)$ associating 
with time and position a celerity is called a \index{traffic 
regulator} \emph{traffic regulator}, feeding traffic celerities 
to mobile's velocities.
\end{Definition}

We thus reserve the word \index{velocity} \emph{velocity} 
$\xi'(t)$ to the   mobiles  which   governs\footnote{If the 
cinematics of the mobile is described by a controlled 
differential equation $\xi'(t)=g_{\xi}(t,\xi(t),u(t))$, the 
microscopic feedback of the mobile is a solution to the equation 
$g_{\xi}(t,p,u)=r(t,p)$ equating velocities and celerities.} 
their evolution by  the cinematic law 
\begin{equation} \label{e:}
\forall \; t , \; \; \xi'(t)\; = \; r(t,\xi(t))
\end{equation}
and the word \index{celerity} \emph{celerity} $r(t,p)$ as an 
information   \emph{attached to the position $p$ at time $t$} of 
the network, \emph{independent of the actual mobile  passing 
through $\xi(t)=p$ at position $p$ at time $t$} using this 
information provided by the traffic regulator.  

For instance, if the mobile is governed by its own dynamical 
second-order differential equation 
$x''(t)=g_{\xi}(t,\xi(t),\xi'(t))$, then its actual evolution on 
the network is governed by the second-order equation 
\begin{equation} \label{e:}   
\forall \; t, \; \; \xi''(t) \; = \; g_{\xi}(t, 
\xi(t),r(t,\xi(t))) 
\end{equation}
for adapting its evolution to the route of the network on which 
it travels.

\subsection*{Staying on the Route}

The first \emph{purpose} is for any mobile evolving on the route 
$K \subset  \mathbb{R}^{p} $ of a network is to be 
viable\footnote{This ``viability constraints'' is not a 
constraint whenever $K := \mathbb{R}$ is a straight road in 
$\mathbb{R}$.}  in the   sense that 
\begin{equation} \label{e:}   
\forall \; t, \; \; \xi (t) \; \in  \; K 
\end{equation}
 
Viability theorems provide an answer whenever the velocities  
$\xi'(t)$ of the mobiles are bounded by an \emph{a priori} limit 
$c(t,\xi(t))$.  Denoting by $T_{K}(p)$   the ``tangent cone'' to 
$K$ at $p \in K$, the viability theorem  requires that the 
\emph{cinematic version} of this constraint is satisfied: 

\begin{equation} \label{e:}   
\forall \;  t  , \; \; \xi'(t) \; \in \;  T_{K}(\xi(t))
\end{equation} 
Therefore, we have to require that the celerity advised to monads 
(see Definition~\ref{d:Monade}, p.\pageref{d:Monade}) is 
constrained by 
\begin{equation} \label{e:}   
\forall \; t \in \mathbb{R}, \; \forall \; p \in K, \;\; r(t,p) 
\leq \; c(t,p) \;\mbox{\rm and } \; r(t,p)  \in \; T_{K}(p)
\end{equation}
Many other requirements, besides this minimal one, could be 
added, providing a more and more restrictive list of 
requirements.

The question, ``design a traffic regulator'', asked, it remains 
to answer it, for given purposes.

\subsection*{Incoming and Outgoing Durations}

Chronological time $t \in \mathbb{R}^{}$ is not sufficient for 
studying the evolution of mobiles, the concept of travel duration 
is mandatory  for ``chaperoning''\footnote{as proteins chaperon 
other proteins in biology, or old aunts their young and innocent 
nieces.\color{Black}} their evolution. Duration function (with 
variable velocities were introduced in 
\cite[Aubin]{aub-12-Durance} and \emph{Time and Money. How Long 
and How Much Money is Needed to Regulate  a Viable  
Economy},\cite[Aubin]{TM} (in economics), but only for 
``outgoing'' durations.

The idea is   to introduce non negative fonctions from $t \mapsto 
d(t) \in \mathbb{R}^{}_{+} $  with variable velocities, denoting 
the \emph{duration needed to reach a position or after starting 
from a position}. For simplifying the exposition, we assume that 
the durations' accelerations are equal to $0$, so that their 
fluidities are constant. 

\begin{Definition} 
\symbol{91}\textbf{Travel Duration with Constant 
Fluidities}\symbol{93}\label{}\index{travel duration} Travel 
durations with constant fluidities are fonctions of the type 
\begin{equation} \label{e:}   
d_{a\varphi}:  t \mapsto  d_{a\varphi}(t)\; := \; \max(0,   a 
(\varphi t-D)) \;\mbox{\rm where}\; a \in 
 \{-1,0,+1\} \;\mbox{\rm and the fluidity}\; \varphi > 0 
\end{equation}
and where $\Omega := \frac{D}{\varphi}$ is the \index{aperture} 
\emph{aperture}   associated with the duration  
$d_{a\varphi}(\cdot)$.\\    We shall say that a travel duration 
is \begin{enumerate}   \item  \index{incoming duration} 
\index{incoming duration}  \emph{incoming} if $a=-1$, i.e.,  if 
the duration is decreasing and reach $0$ at $\Omega := 
\frac{D}{\varphi}$ for $\varphi >0$;

 \item   \index{stationary} \emph{stationary}  if $a=0$;  

\item \emph{outgoing} if $a=+1$, i.e.,  if the duration is 
increasing and starts from  $0$ at $\Omega := \frac{D}{\varphi}$ 
for $\varphi >0$.
\end{enumerate}  
\end{Definition}

One motivation for introduce these durations is that they signal 
\index{junction time} \emph{junction times} $\Omega $ when 
$d(\Omega )=0$ and that this information is provided to dynamical 
system when, for instance, a differential inclusion $x'(t) \in 
F(t,x(t))$ involves the duration. 

For example, using constant fluidities $\varphi >0$, we obtain

\begin{enumerate}   
\item  \textbf{Incoming evolutions} governed by

\begin{equation} \label{e:}   
x'(t) \; \in \; \max(0, \Omega  - \varphi t) F(t,x(t))
\end{equation}
which decreases the velocities  of $x(\cdot)$ for $t < \Omega $ 
until they vanish when $t = \Omega $;

\item  \textbf{Outgoing evolutions} governed by
 
\begin{equation} \label{e:}   
x'(t) \; \in \;  \max(0,\varphi t -\Omega ) F(t,x(t))
\end{equation}
which increases the velocities  of $x(\cdot)$ from $0$ when $t = 
\Omega  $ to positive ones when $t > \Omega  $. 
  \end{enumerate} 

These are examples of differential inclusions of the form

\begin{equation} \label{e:}   
x'(t) \in F(t, d(t), x(t))
\end{equation}

General travel durations are non negative functions $t \mapsto 
d(t) \in \mathbb{R}^{}_{+}$ which   successively go through an 

\begin{enumerate}   
\item  \textbf{Incoming phase}, when durations are  decreasing 
from their local maximal to their local minima at \index{junction 
date} \emph{junction date} $\Omega $ where $d(\Omega )=0$ 
required to be equal to $0$;

\item  \textbf{Stationary phase}, when durations remain equal to 
$0$;

 \item \textbf{Outgoing phase} when durations are  
increasing from their local minima $d(\Omega )=0$  to their local 
maxima at \index{reversal date} \emph{reversal date} after which 
they increase.
  \end{enumerate} 

A \index{travel schedule} \emph{travel schedule} is the 
concatenation of decreasing, stationary and increasing durations 
in this order.  

Such situations could be qualified as \emph{``structured by 
durations''} since McKendrick opened the way to equations 
\emph{``structured by ages''}\footnote{The age-structured 
standard approach starts with the establishment of the 1926 
McKendrick partial differential equation relating the population 
and its partial derivatives with respect to time and age. 
\glossary{McKendrick (Anderson Gray ) [1876-1943]} Time and age 
are scalar durations evolving with constant fluidities equal to 
$1$.   Age-structured partial differential equations involving 
both time and age have been extensively studied (see, among an 
extensive literature, \cite[Anita]{Anita}, 
\cite[Aubin]{RegulationBirths-11}, \cite[Aubin, Bonneuil \& 
Maurin]{abb99evs}, \cite[Iannelli]{im95aspd},  \cite[Keyfitz N. 
\& Keyfitz B.]{ke297dem}, \cite[Mckendrick]{mck26dem}, \cite[Von 
Foerster]{foe59dem}, \cite[Webb]{web85age}, etc.).} when age 
plays the rôle of an  outgoing duration   with constant fluidity 
equal to $+1$, and thus,  increasing ... forever (forgetting 
\glossary{Kafka (Franz ) [1883-1924]} \emph{Franz Kafka}'s 
observation that   \emph{``eternity is long, above all towards 
the end''}).

\vspace{ 5 mm} 

\textbf{Remark: Junction detectors} --- Actually, durations can 
be regarded as ``temporal junctions detectors'' for 
characterizing junctions. One can, for instance, also regard then 
as ``spatial junctions detectors'' by defining the function 
$d(\cdot)$ through a distance function $x \mapsto 
\mathbf{d}_{\mathbb{J}}(x)$ to a the junction:

\begin{equation} \label{e:}   
\forall \; t, \; \; d(t) \; := \;  \mathbf{d}_{\mathbb{J}}(p(t))
\end{equation}
Indeed, at time $\Omega $ when $d(\Omega )=0$, then 
$\mathbf{d}_{\mathbb{J}}(p(\Omega ))=0$ so that  $p(\Omega ) \; 
\in \; \mathbb{J} $.

Other junction detectors motivated by such and such concrete 
problem can be devised, the modifications from the temporal 
junction detector being marginal, since this study assumes that 
the evolution of monades are governed by differential inclusions 
$x'(t) \in F(t,d(t),p(t),x(t))$ where $d(\cdot)$ can be any 
  detector, feeding back on any component of the traffic 
state $(t,d,p,x)$. For instance, for spatial junction  detectors, 
$x'(t) \in F(t,d_{\mathbb{J}}(p(t)),p(t),x(t))$. \hspace{ 2 mm} 
\hfill $\;\; \blacksquare$ \vspace{ 5 mm}

\subsection*{From Mobiles to Monades}

In order to construct traffic regulators providing celerities to 
the mobiles passing through a position at a given time, we can 
take into account the knowledge of  some complementary 
\index{traffic attributes}  traffic attributes $x$ attached to 
the routes of the networks, \index{monad} that we shall call
\emph{monads}\footnote{From the Greek \emph{``monos''},  
\emph{unique}. \glossary{Leibniz (Gottfried) [1646-1716]} 
\emph{Gottfried Leibniz}  already used the concept of 
\emph{monad} in philosophy, and the name was borrowed in  
mathematical category theory and abstract programming}  to place 
under a same abstract umbrella those many examples of attributes 
of traffic information provided to ``mobiles'',  the physical 
nature of which is not involved  besides the fact that they 
circulate along the routes of the  network. 

For instance, monads can represent the celerity or the measure of 
jam at each time $t$ and each position $p$. When we want to 
describe a position $p$ at which the mobiles must stop at time 
$t$, we have to impose that the celerity (monad) is equal to $0$ 
at this time and position. Also, we may define several measures 
of \emph{jam} before or behind a given position, used as another 
monad, and so on.

\begin{Definition} 
\symbol{91}\textbf{Monades}\symbol{93}\label{d:Monade}\index{monad} 
A network $K \subset \mathbb{R}^{p}$, being given, we  add to the 
triplets $(t,d,p) \in  \mathbb{R}^{}\times \mathbb{R}^{}_{+} 
\times  K $ describing a chronological time, a duration and a 
position other traffic attributes $x \in X:=\mathbb{R}^{m}$, 
regarded as the \index{monad map} \emph{monad space}, that we 
shall call \index{monads} \emph{monads}, for completing the 
description of a \index{traffic state} \emph{traffic state}
 $(t,d,p,x) \in  \mathbb{R}^{}\times \mathbb{R}^{}_{+} 
\times  K \times \mathbb{R}^{m}$.  A \index{monad map} 
\emph{monad map} $\mathbb{M}: (t,d,p) \in \mathbb{R}^{}\times 
\mathbb{R}^{}_{+} \times  K  \leadsto  \mathbb{M}(t,d, p) \subset 
X$ associates with each triplet $(t,d,p)$  a subset 
$\mathbb{M}(t,d, p)$ of monads $x$. A \index{monad regulator} 
\emph{monad regulator} is a set-valued map associating with each 
traffic state $(t,d,p,x) \in \mbox{\rm Graph}(\mathbb{M})$ a 
subset $R(t,d,p,x) \subset  \mathbb{R}^{m}$ of \emph{celerities} 
advised to the mobiles 
passing through the traffic state $(t,d,p,x)$.\\
The derivative $t \mapsto x'(t)$ of the evolution of monads is 
regarded as the \index{surge} \emph{surge} of the traffic at time 
$t$. 
\end{Definition}

The monads evolve   along the routes of the network according to 
some evolutionary laws, differential equations or inclusions for 
accommodating controls and/or tyches (disturbances, 
perturbations).   Traffic regulators provide the information on 
celerity to mobiles between two junctions  and the discontinuous 
behavior at junctions is described by \index{junction regulator} 
\emph{junction regulators} described later, after we describe two 
examples of traffic attributes. For instance,

\begin{enumerate}

\item   \textbf{\emph{Dynamic Traffic Regulator}} 

The simplest example of monad is the celerity  $x(t)=p'(t)$ 
advised to the mobiles along the routes. Hence   the ``surge'' 
boils down to the acceleration $x'(t)=p''(t)$. Introducing these 
monads is mandatory whenever constraints on the velocities of the 
mobiles are used, for instance whenever mobiles have to stop at 
given times and positions of the route\footnote{For instance, in 
road networks, stop signs and traffic lights requiring to stop 
during periodic periods of time.}.   For instance, if $p=1$ and 
if  the surge $x'(t):=p''(t)$ is equal to $f(t,d,p,r)$, the 
traffic regulator $r(t,p)$   at time $t$ at position $p$ 
providing the celerity is a solution to the Burgers type partial 
differential equation 
\begin{equation} \label{e:}   
\frac{\partial r(t,d,p)}{\partial t}+ \varphi \frac{\partial 
r(t,d,p)}{\partial d} +r(t,d,p) \frac{\partial r(t,d,p)}{\partial 
p} \; = \; f(t,d,p,r(t,d,p))
\end{equation} 
associated with adequate Cauchy and/or Dirichlet conditions (see, 
for instance, \cite[Aubin]{AUB-10-STR} and Chapter~16, p. 631, of 
\emph{Viability Theory.  New Directions}, \cite[Aubin, Bayen \& 
Saint-Pierre]{absp} and, in an economic framework, Chapter~2, p. 
31, of \emph{Time and Money. How Long and How Much Money is 
Needed to Regulate  a Viable  Economy},\cite[Aubin]{TM}.

\item  \textbf{\emph{Measures of Traffic Jam}} 

In another instance, we may require  that the  jam 
(«~engorgement~» in French) of the traffic satisfies some 
  constraints. 

By definition, the \index{traffic jam} (measure of) \emph{traffic 
jam} at position $p \in K$ is  \emph{the number $x(p)$ of mobiles on a 
stretch of route before   and/or after 
 position} \footnote{described for instance  by $K \cap \gamma B \cap 
(x\pm \mathbb{R}^{\ell}_{+})$. Instead of the number $x(p)$ of 
mobiles in this set, we can also compute the measure occupied by 
the mobiles in this stretch of road.} $p \in K$, time $t$ and 
duration $d$.

The derivative $x'(t) \in X^{\star}$ of a differentiable jam 
function has the dimension of a \index{surge} \emph{surge}.  We 
can require for instance that at each time $t$, the \index{jam 
function} \emph{jam function} or  \emph{bottleneck  function} 
$(t,d,p) \in  b (t,d , p ) \subset \mathbb{B}(t,d,p) \subset  
\mathbb{R}^{m}$  is a ``jam capacity'' that the route allows at 
time $t$ for a duration $d$  at position $p$. When $p=1$, we 
recover the characteristic system underlying the 
Hamilton-Jacobi-Moskowitz summarized in Chapter~14 of 
\emph{Viability Theory.  New Directions}, \cite[Aubin, Bayen \& 
Saint-Pierre]{absp}, and studied in an abundant literature. 
\hfill $\;\; \blacksquare$ \vspace{ 5 mm}
 
 \end{enumerate}
 
\emph{These requirements are part of a list that traffic modelers 
provide, from which we must built a traffic regulator providing 
mobiles their celerity at all positions and instants, information 
which they should or must use for satisfying them. }

\subsection*{Junctions}

The evolutions of mobiles on a  network  between incoming   
 and outgoing states  can be  
\emph{``punctuated''} at junctions of the network, where  
\index{intermodal systems} \emph{``intermodal systems''} between 
two prejunction and postjunction states of a  \index{junction 
relation} \emph{junction relation}. For instance,a network of 
routes is a set of routes punctuated by junctions, in particular, 
by impulse junctions such as cross-roads, synapses, clearing 
houses, bit nodes, etc. The circulation between two departure and 
arrival states on a network  can be interrupted by an 
\emph{intermodal system} between \emph{prejunction} and 
\emph{postjunction} dates and positions, whereas at crossings, 
these prejunction and postjunction dates collapse at a same 
\emph{impulsive junction} date and position. 

In the case, of biological networks,  the transmittal of proteins 
and other chemical compounds (hormones, etc.)  are produced in 
endocrine glands for regulating specific physiological processes  
by travelling  through the bloodstream    between two distant 
emitter  and receptor. They   last between two prejunction and 
postjunction dates whereas, in neural networks, for instance, 
other proteins (neurotransmitters) crossing a synapse between a 
presynaptic neuron and a postsynaptic one trigger 
``instantaneously'' the transportation of neurotransmitters at an 
impulsive (synaptic) junction  date.

Regarding the monads, even though the position of the junction 
may be the same, the monads   ``passing'' from one attribute to 
another one  during a junction, \index{impulse}  impulsive (when 
it is triggered by an   \emph{impulse}) or not.  Therefore, the 
description of a network of routes is not restricted to the 
continuous time circulation of monads between outside  junctions, 
but integrate also   discontinuous evolutions at junctions, 
governed by   junction maps   operating at junctions, impulsive 
or not. 

For instance, if the traffic attribute of a monad is the 
celerity, it may be equal to $0$ for some monads on some routes 
at some time to stop at the junction whereas the monads on other 
routes have positive celerities at some time ans switch to a 
different configuration at an other time (in road traffic, for 
instance). If the traffic attribute is a (measure of) jam, at 
some junction and at some time,   the jam of some monads have to 
switch (immediately, in the impulsive case) for the sum of the 
jams to be below a given threshold. 

\emph{How can we locate or signal a junction state?} The 
introduction of duration functions allows the system to detect 
then whenever the durations vanish. 

Hence, as for impulse control and hybrid systems, the behavior of 
the monads at the junction is described by a \index{junction 
relation} \emph{junction relation} , either single-valued or 
set-valued to accommodate controls and/or tyches (perturbations, 
disturbances, etc.): see  Definition~\ref{d:JunctionRelation}, 
p.\pageref{d:JunctionRelation} of Section~\ref{s:junction}, p. 
\pageref{s:junction}.

\subsection*{Motivation: Road Traffic}

Road traffic information exists at least since the milestones 
have been used by Roman Empire for deducing the travel durations 
  between two positions\footnote{In his 1960's study of 
the \emph{«~yamonamö~»}, still a stone age tribes at this time,   
\emph{Napoleon Chagnon} reports in \cite[\emph{Noble 
Savages}]{Chagnon} that the distance between two places is 
measured by durations (number of ''sleeps'' during a trip). The 
yanomamö cannot count accurately beyond $2$ and had  to sit down 
to estimate a long distance by using their toes. With milestones 
and the concept of velocity, travel duration could be estimated 
from the information on distances, a spatial metaphor of time, 
actually, duration.}. At the time, only the position was 
indicated, the celerity being not yet an issue! Nowadays, road 
traffic regulation, as rudimentary and frustrating it is, is 
known to everyone: they are provided by speed limit  signs (an 
upper bound of the velocity used by the mobiles passing by), or 
more coercitive information, as traffic lights or signals, 
imposing a speed limit equal to $0$ during a given 
time-interval,  or other signs providing qualitative advices, 
etc. \emph{Such information is provided by a road equipment to 
the vehicles\footnote{The recommended celerities   can be posted  
on VMS (variable message signs),   broadcasted on mobile phones 
equipped with GPS,  displayed on  twinned speedometers both the 
effective velocity and the recommended celerity   
 (their difference triggering alarms), regulated by 
cruise control systems  adjusting automatically the velocity to 
be equal to the broadcast  celerity, etc.}, which have to abide 
by this information for moving on the road.}


The construction of  feedbacks piloting evolutions viable in an 
environment, here, a network of routes,  is the very objective of 
viability theory, the results of which allow us to compute the 
celerities regulating the traffic of monads in this network.   

\mbox{}

\emph{The search of celerities is a by-product of most studies of 
road traffic}. \glossary{Greenshields (Bruce D.) []} 
\emph{Greenshields} used in 1933 photographic measurement methods 
for the first time to describe a \emph{phenomenological law} 
described by  a quadratic relation between vehicles and their 
density and flows, and using it as inputs of first-order partial 
differential equations (conservation laws and 
Hamilton-Jacobi-Bellman equations) providing ``cumulative 
number''   of vehicles passing at given  position after a given 
time   on a one-dimensional  road (see \cite[Aubin, Bayen \& 
Saint-Pierre]{absp08hj},
 \cite[Aubin]{AUB-Leit09,Lax-Hopf}, Chapter~14, p. 563, of  \emph{Viability Theory.  New Directions}, \cite[Aubin, Bayen \& 
Saint-Pierre]{absp}, devoted to regulation of traffic and the 
forthcoming \cite[\emph{Mathematical Approaches to Traffic 
Management}]{ADA-Traffic}). Furthermore, the viability algorithms 
allow us to compute the traffic function and to regulate viable 
evolutions (see \cite[Bayen, Claudel, Saint-Pierre]{BCSP07b, 
BCSP07a}, \cite[Claudel, Bayen]{CB08a,CB08b}), 
\cite[Désilles]{Anya}. 

This \emph{``shift from densities to celerities''} was hidden, so 
to speak, in  the description  of the traffic by the number of 
vehicles. They are related by a function (called the fundamental 
diagram) relating densities to flows by a system  of  partial 
differential equations. The Legendre-Fenchel transform of density 
function is a function associating flows with celerities, 
regarded as controls of a dynamical system governing the 
evolutions of positions.  \glossary{Daganzo (Carlos F.) []} 
\emph{Daganzo} discovered celerities (see 
\cite[Daganzo]{DaganzoCell94,DaganzoCell95,DaganzoVP-1,DaganzoVP-2}) 
in his ``Daganzo variational principle" without naming them for 
defining optimal evolutions.

Right or wrong, we choose the other way around, \emph{finding 
traffic regulators providing celerities first, and, in the 
classically formulated problems,   uncovering as a by-product the 
underlying first-order partial differential equations}, just to 
provide a link with the eighty years abundant traffic engineering 
literature, but not for using them. Because partial differential 
equations do not provide more information on the regulation  than 
the one directly provided by the viability approach regarding the 
mathematical, algorithmical and software issues.

\subsection*{Contents}

Section~\ref{s:Relations}, p. \pageref{s:Relations} defines 
\emph{relations}. Usually a set of variables is decomposed into 
input variables and output variables, and (set-valued) maps sends 
inputs into outputs. However, set-valued analysis favors a 
graphical approach regarding maps as their graphs in the product 
of the input  and the output spaces. This dichotomy between 
inputs and outputs is not always reasonable: what matters is a 
subset of variables ranging over a product of two  \emph{or more} 
 spaces, playing the rôle of the graph in the case of two 
spaces. Hence all results on graphs of maps can be adapted, 
except, naturally,  the ones which are formulated in terms of 
set-valued maps. Since there are many variables describing a 
traffic state, and twice as many when we consider junctions, it 
seamed necessary to abandon set-valued maps in favor of relations 
on two or more spaces. Consequently, in Section~\ref{s:junction}, 
p. \pageref{s:junction}, we define transport relations for 
singling out incoming and outgoing states described by time, 
duration, position and monades, i.e., eight variables, and, even 
more, if the monades themselves are split in relevant classes of 
variable (celerities, jam, etc.). Hence transport relations are 
subsets of the product of eight relevant spaces and junctions 
relations are transport relations with vanishing durations. 
Section~\ref{s:TransportEvolutions}, p. 
\pageref{s:TransportEvolutions}  defines transport evolutions 
linking an incoming state to an outgoing one across a junction. 
Section~\ref{s:TransportKernelRelation}, p. 
\pageref{s:TransportKernelRelation} adresses the problem under 
investigation: \emph{how can we construct a transport relation 
across a junction linking incoming and outgoing states?} This is 
kind of \emph{geodesic} problem, the formulation of which is 
formulates as the \index{transport kernel relation} 
\emph{transport kernel relation}, i.e.,  the subset of 
``linkable'' pairs of incoming and outgoing states across the 
junction. Section~\ref{s:ConstructionTransKern}, p. 
\pageref{s:ConstructionTransKern} characterizes transport kernel 
relations as capture basins of an auxiliary system. Therefore, 
transport kernel inherits the properties of capture basins. Among 
them, the   tangential conditions characterizing capture basins 
provide the transport regulators we were looking for, as it is 
explained in Section~\ref{s:TransportRegulator}, p. 
\pageref{s:TransportRegulator}.

\section{Relations} \label{s:Relations}

Since traffic states involve   time, duration, position and 
monad  components, what matters is the subset of such states in 
the product $\mathbb{R}^{}\times \mathbb{R}^{}_{+}\times 
\mathbb{R}^{p} \times \mathbb{R}^{m}$. Until now, we regarded 
such a subset as the graph $\mbox{\rm Graph}(\mathbf{M})$ of the 
monad map  $\mathbb{M}: \mathbb{R}^{}\times 
\mathbb{R}^{}_{+}\times \mathbb{R}^{p}  \leadsto  
\mathbb{R}^{m}$. Although reasonable, this choice is arbitrary, 
and we could as well, by ``reorganizing'' this product as 
$\mathbb{R}^{}\times \mathbb{R}^{p} \times \mathbb{R}^{m}\times 
\mathbb{R}^{}_{+}$, regard the same subset as the graph of the 
duration $\mathbb{D}: \mathbb{R}^{}  \times \mathbb{R}^{p} 
\times  \mathbb{R}^{m} \leadsto \mathbb{R}^{}_{+}$ associating 
with triplet time, position and monad its duration. 

Moreover, transport  relation  deals  with two traffic states, 
incoming state $(t_{in},d_{in},p_{in},x_{in})$ and  outgoing 
state $(t_{ou},d_{ou},p_{ou},x_{ou})$,  constituting a 
\index{transport state} \emph{transport state}

\begin{equation} \label{e:}   
((t_{in},d_{in},p_{in},x_{in}),(t_{ou},d_{ou},p_{ou},x_{ou})) \; 
\in \; (\mathbb{R}^{}\times \mathbb{R}^{}_{+}\times 
\mathbb{R}^{p} \times \mathbb{R}^{m}) \times (\mathbb{R}^{}\times 
\mathbb{R}^{}_{+}\times \mathbb{R}^{p} \times \mathbb{R}^{m})
\end{equation}

Consequently, it seems pointless to privilege one decomposition 
of this product as the product of two partial products of theses 
spaces, since there are too many combinations separating what are 
the input variables and output variables.

The incoming-outgoing ``natural'' breeching would be to regard a 
subset of  $\mathcal{\mathcal{Q}} \subset (\mathbb{R}^{}\times 
\mathbb{R}^{}_{+}\times \mathbb{R}^{p} \times \mathbb{R}^{m}) 
\times (\mathbb{R}^{}\times \mathbb{R}^{}_{+}\times 
\mathbb{R}^{p} \times \mathbb{R}^{m})$, not only as the graph of 
the map

\begin{equation} \label{e:}   
((t_{in},d_{in},p_{in},x_{in})  \leadsto 
(t_{ou},d_{ou},p_{ou},x_{ou})) 
\end{equation}
but also as the graph of the map
\begin{equation} \label{e:}   
(((t_{in},d_{in},p_{in}),x_{in}),(t_{ou},d_{ou},p_{ou}))  
\leadsto 
 x_{ou} 
\end{equation}
providing the arrival monad. Set-valued analysis and viability 
theory use systematically the graphical approach for dealing with 
set-valued maps (or the epigraphical approach for studying 
extended numerical functions). An abundant literature concerns 
this view point (see among many monographs \emph{Variational 
Analysis}, \cite[Rockafellar  \& Wets]{rw91nsa}, \emph{Set-valued 
analysis} , \cite[Aubin \& Frankowska]{af90sva}, \emph{Viability 
Theory.  New Directions}, \cite[Aubin, Bayen \& 
Saint-Pierre]{absp}, etc.).  

This is the reason why we shall define \index{relation} 
\emph{relations}, generalizing  graphs of set-valued maps, 
regarded as \emph{binary relations}. The idea of maps is so 
entrenched in our minds that it seems difficult to go from the 
number two\footnote{At this stage, mathematics are in the same 
situation than the yamonamö who have a name for the number two 
only, as we have a name only for maps, which are binary 
relations.} to higher numbers of elements, without separating 
variables in two categories only, inputs and outputs, top and 
down, left or right.

\begin{Definition} 
\symbol{91}\textbf{Relations}\symbol{93}\label{d:Relations} 
\index{relations} A relation $\mathcal{R} \subset  
\prod_{i=1}^{\ell} X_{i}$ is a subset of the product of $\ell $ 
spaces. It relates the multistates $(x_{i})_{i=1, \ldots ,\ell}$ 
belonging to $\mathcal{R}$.  A relation $\mathcal{R}_{1}$ is a 
restriction of a relation $\mathcal{R}_{2}$ (and 
$\mathcal{R}_{2}$ is an extension of $\mathcal{R}_{1}$) if   
$\mathcal{R}_{1}$ is contained $\mathcal{R} \subset \prod_{i \in \mathbb{I}}^{}X_{i}$in the graph of 
$\mathcal{R}_{2}$.\\
The \index{hyperrelation} \emph{hyperrelation} $ 
\overbrace{\mathcal{R}}$ of a relation $\mathcal{R}$ is a 
relation  $\overbrace{\mathcal{R}} \subset  
\prod_{i=1}^{\ell}\mathcal{P}(X_{i}) $ made of families 
$(A_{i})_{i=1, \ldots ,\ell} $ such that $\prod_{i 
=1}^{\ell}A_{i} \; \subset \; \mathcal{R}$.  
\end{Definition}
 
 
\section{Transport  and Junction 
Relations} \label{s:junction}

A \emph{transport relation} involves two monads $\displaystyle{ 
x_{in}  \in  \mathbb{M} (T_{in},D_{in} , P_{in} )}$ at 
\emph{``departure''} date $T_{in}$ before the 
\emph{``prejunction''} date $\Sigma_{in}$ of the junction and 
$\displaystyle{ x_{ou}(t) \in  \mathbb{M} (T_{ou},D_{ou} , 
P_{ou})}$ at \emph{``arrival''} date $T_{ou}$ after the 
\emph{``postjunction''} date $\Sigma_{ou}$ of the junction.  

We denote by $\mathbb{R}^{2}_{\leq }$ the subset of pairs 
$(T_{in},T_{ou})$ of departure and arrival dates (or pairs 
$(\Sigma_{in} ,\Sigma_{ou})$ of prejunction and postjunction 
dates)  such that $ T_{in} \leq T_{ou}$ (or $ \Sigma_{in} \leq 
\Sigma_{ou}$). We shall naturally \emph{assume once and for all 
that the departure date $T_{in} \leq T_{ou}$ is always smaller 
than the arrival date $T_{ou}$ as well as  the prejunction date 
$\Sigma_{in} \leq \Sigma_{ou}$ is always smaller than the 
postjunction date}.

The incoming decreasing duration $d_{in}(t)$ and outgoing 
increasing duration $d_{ou}(t)$ are used to \emph{locate 
junctions} whenever $d_{in}(\Sigma_{in})=0$ at the prejunction of 
the junction and  $d_{ou}(\Sigma_{ou})=0$ at the postjunction of 
the junction. 

The ``continuous time traffic''   of monads along the route 
$\mathbb{M}:(t,d,p)  \leadsto \mathbb{M}(t,d,p)$ is interrupted 
(punctuated) and replaced by an \index{intermodal transit} 
\emph{``intermodal''}  transit governed by a given different 
intermodal system. 

\begin{Definition} 
\symbol{91}\textbf{Transport     
Relations}\symbol{93}\label{d:junctionPunctuationMap} A 
\index{transport relation} \emph{transport relation} is a 
relation $\mathcal{Q}$ contained in the product
\begin{equation} \label{e:}   
\mathcal{Q} \; \subset  \; (\mathbb{R}^{}\times 
\mathbb{R}^{}_{+}\times \mathbb{R}^{p} \times \mathbb{R}^{m}) 
\times (\mathbb{R}^{}\times \mathbb{R}^{}_{+}\times 
\mathbb{R}^{p} \times \mathbb{R}^{m})
\end{equation}
which relates incoming and outgoing states
\begin{equation} \label{e:}   
((t_{in},d_{in},p_{in},x_{in}),(t_{ou},d_{ou},p_{ou},x_{ou})) \; 
\in \; \mathcal{Q}
\end{equation}
For instance, the \index{} \emph{monad relation} is defined by 
\begin{equation} \mathcal{M} \; := \;  \mbox{\rm 
Graph}(\mathbb{M})^{2} \end{equation} We shall say that a 
\emph{transport relation is viable in the monad relation} of a 
network if $\mathcal{Q}  \subset \mathcal{M}  $ is contained in 
the monad relation $\mathcal{M}$, product of the graphs of the 
monad maps.
\end{Definition}

Junction  relations are transport  relations with vanishing 
incoming and outgoing  durations, which allow  them to be 
detected by incoming and outgoing durations functions:

\begin{Definition} 
\symbol{91}\textbf{Junction  
Relations}\symbol{93}\label{d:JunctionRelation}\index{junction  
relation} A \emph{junction  relation} $\mathcal{J}$ is defined as 
transport relation satisfying

\begin{equation} \label{e:JunctionRelation}   
\mathcal{J} \; \subset  \; (\mathbb{R}^{}\times \{0\}\times 
\mathbb{R}^{p} \times \mathbb{R}^{m}) \times (\mathbb{R}^{}\times 
\{0\}\times \mathbb{R}^{p} \times \mathbb{R}^{m})
\end{equation}
when the incoming and outgoing durations of which are equal to 
$0$   and when
\begin{equation} \label{e:} \left\{ \begin{array}{l}  
\Sigma^{\mathcal{J}}_{in}  \; := \; 
\sup_{((\Sigma_{in},\Pi_{in},\Xi_{in}),(\Sigma_{ou}, 
\Pi_{ou},\Xi_{ou})) \in \mathcal{J}}\Sigma_{in} \\
\leq \\ \Sigma^{\mathcal{J}}_{ou}  := 
\inf_{((\Sigma_{in},\Pi_{in},\Xi_{in}),(\Sigma_{ou}, 
\Pi_{ou},\Xi_{ou})) \in \mathcal{J}}\Sigma_{ou} 
\end{array} \right. 
\end{equation} \\
The triplet $(\Sigma_{in}, \Pi_{in},\Xi_{in}) \in \mathbb{R}^{} 
\times \mathbb{R}^{p} \times \mathbb{R}^{m}$ is called 
\index{prejunction} \emph{prejunction}  state at date 
$\Sigma_{in} \in \mathbb{R}^{}$ and position  $\Pi_{in}$ and the 
triplet $(\Sigma_{ou}, \Pi_{ou},\Xi_{ou}) \in \mathbb{R}^{} 
\times \mathbb{R}^{p} \times \mathbb{R}^{m}$  is a
\index{postjunction} \emph{postjunction} state at date  
$\Sigma_{ou} \geq \Sigma_{in}$ date and position $\Pi_{ou}$.\\
A junction is called  \index{time impulsive junction} \emph{time 
impulsive} if the prejunction and postjunction   dates  
coincide:  $ \Sigma_{in}= \Sigma_{ou} =: \Sigma$.  
\end{Definition}

We shall identify prejunction and postjunction  junctions states 
 \begin{equation} \label{e:}   
((\Sigma_{in},0,\Pi_{in},\Xi_{in}),(\Sigma_{ou},0,\Pi_{ou},\Xi_{ou})) 
\; \in \; \mathcal{J}
\end{equation} 
with a relation  
 \begin{equation} \label{e:}   
((\Sigma_{in}, \Pi_{in},\Xi_{in}),(\Sigma_{ou}, 
\Pi_{ou},\Xi_{ou}))  \; \in \; \mathcal{J} \subset (\mathbb{R}^{} 
\times \mathbb{R}^{p} \times \mathbb{R}^{m}) \times 
(\mathbb{R}^{} \times \mathbb{R}^{p} \times \mathbb{R}^{m})
\end{equation}  

The state of a incoming-outgoing evolution at a junction satisfies

\begin{equation} \label{e:}   
(\Sigma_{in},0, 
p_{in}(\Sigma_{in}),x_{in}((\Sigma_{in})),(\Sigma_{ou},0, 
p_{ou}(\Sigma_{ou}),x_{ou}(\Sigma_{ou})))\; \in \; \mbox{\rm 
Graph}(\mathcal{J})  
\end{equation}

Therefore, at a junction relation, incoming and outgoing monad 
evolutions satisfy

\begin{equation} \label{e:} \left\{ \begin{array}{l}  
x_{in}(\Sigma_{in}) \; \in \;  \mathbb{M}(\Sigma_{in},0, 
p_{in}(\Sigma_{in}))\\
x_{ou}(\Sigma_{ou}) \; \in \;  \mathbb{M} (\Sigma_{ou},0, 
p_{ou}(\Sigma_{ou}))\\
\end{array} \right. 
\end{equation}
or, equivalently. 
\begin{equation} \label{e:}   
(\Sigma_{in},0, 
p_{in}(\Sigma_{in}),x_{in}((\Sigma_{in})),(\Sigma_{ou},0, 
p_{ou}(\Sigma_{ou}),x_{ou}(\Sigma_{ou})))\; \in \; \mbox{\rm 
Graph}(\mathbb{M}) ^{2}
\end{equation}

The simple junctions are the junctions $\mathcal{J}= 
\mathcal{J}_{in}\times \mathcal{J}_{ou}$ where

\begin{equation} \label{e:} \left\{ \begin{array}{l} 
\mathcal{J}_{in} \subset \mbox{\rm Graph}(\mathbb{M}) \cap 
(\mathbb{R}^{}\times \{0\}\times 
\mathbb{R}^{p} \times \mathbb{R}^{m}) \\
\mathcal{J}_{ou} \subset \mbox{\rm Graph}(\mathbb{M}) \cap 
(\mathbb{R}^{}\times \{0\}\times \mathbb{R}^{p} \times 
\mathbb{R}^{m})
\end{array} \right. 
\end{equation}

In this case, the junction crossing is split into two parts: 
arriving at the target $\mathcal{J}_{in}$ and leaving from the 
source $\mathcal{J}_{ou}$, \emph{independently}. The study of 
these two problems is a consequence of the more realistic case 
when the junction relation is not a product of subset, so that 
the prejunction states influence  the postjunction states.

\section{Transport Evolutions Crossing a Junction Relation} 
\label{s:TransportEvolutions}

From now on, we assume that the accelerations of the incoming and 
outgoing durations are equal to $0$, to that their fluidities 
(absolute velocities) $\varphi_{in}> 0$ and $\varphi_{ou} > 0$ 
are constant. These fluidities are parameters of the problem. The 
\index{aperture} \emph{aperture} $\Omega \geq 0$  measures 
\emph{the time spent before reaching and after leaving a 
junction}.

For a given aperture $\Omega $, the transport of monads from the 
departure date $T_{in}$ to the arrival date $T_{ou}$ is 
decomposed into three phases: 

\begin{enumerate}   
\item an \emph{incoming monad evolution} $t \mapsto x_{in}(t)$ 
starting from departure monad $x_{in}$ at departure time $T_{in}$ 
and arriving at the prejunction monad $x_{in}(T_{in}+\Omega)$ at 
the  prejunction date $\Sigma_{in}=T_{in}+\Omega $;
  
\item an \emph{intermodal evolution} during a \index{intermodal 
period} \emph{intermodal period} (or \emph{junction time 
interval} $[\Sigma_{in} ,\Sigma_{ou}]$), which is governed by a 
\emph{``junction relation''};

\item  an \emph{outgoing monad evolution} $t \mapsto x_{ou}(t)$ 
starting at postjunction date $\Sigma_{ou}:=T_{ou}-\Omega$ at the 
postjunction monad $x_{ou}(T_{ou}-\Omega)$ and arriving at time 
$T_{ou}$ at the arrival monad $x_{ou}$. 
\end{enumerate} 

\begin{Definition} 
\symbol{91}\textbf{Transport Evolution Across a 
Junction}\symbol{93}\label{d:TransportEvo}\index{} Let us 
consider the monad relation $\mathcal{M}$ describing a network, a 
viable junction  relation $\mathcal{J} \subset \mathcal{M}$ of 
the network and  incoming and outgoing fluidities $\varphi_{in}$ 
and 
$\varphi_{ou}$. \\
A pair of incoming and outgoing evolutions $t \mapsto 
((t,d_{in}(t),p_{in}(t),x_{in}(t)), 
(t,d_{ou}(t),p_{ou}(t),x_{ou}(t)))$ is a \index{transport 
evolution} \emph{transport evolution starting at 
$(T_{in},D_{in},P_{in},x_{in})$ and arriving at 
$(T_{ou},D_{ou},P_{ou},x_{ou})$ across the junction 
$\mathcal{J}$}  if there exist  

\begin{enumerate}   \item  an 
\index{aperture} \emph{aperture} $\Omega \geq 0$;  \item   a 
junction state $((\Sigma_{in}, \Pi_{in},\Xi_{in}),(\Sigma_{ou}, 
\Pi_{ou},\Xi_{ou})) \in \mathcal{J}$  
\end{enumerate}
satisfying the following requirement:  their 
\index{concatenation}  $\Omega $-\emph{concatenation} evolution 
$t \mapsto x(t):=(x_{in} \diamondsuit_{\Omega } x_{ou})(t)$ of 
the monad between the departure date  $T_{in}$ and the arrival 
date  $T_{ou}$ of the incoming evolution $t \in 
[T_{in},T_{in}+\Omega ] \mapsto  x_{in}(t) $ and the outgoing 
evolution $t \in  [T_{ou}-\Omega ,T_{ou}] \mapsto x_{ou}(t)$ 
satisfies the properties summarized in the tables below: 
\scriptsize 
\begin{equation} \label{e:}   
\begin{array}{c|c|c|c|c|c|c} 
  & \mbox{continuous}& prejunction  & \mbox{junction} &  postjunction &\mbox{continuous} &  \\
    & \mbox{evolution}& state  & \mbox{jump} & state  &\mbox{evolution} &  \\ \hline 
  & t \in [T_{in},T_{in}+\Omega ]&   & \mbox{jump} &   &t \in 
  [T_{ou}-\Omega ,T_{ou}] &  \\       
\hline T_{in} &  t & \Sigma_{in}:=T_{in}+\Omega & \Rrightarrow & 
\Sigma_{ou}= T_{ou}-\Omega &t& T_{ou}\\   D_{in} & D_{in}+ 
\varphi_{in}(T_{in}-t) & 0 & \Rrightarrow & 0 &D_{ou} 
+\varphi_{ou}(t-T_{ou})  & D_{ou}\\ P_{in} & p_{in}(t) & \Pi_{in} 
= p_{in}(T_{in}+\Omega)
 & 
\Rrightarrow & \Pi_{ou}= p_{ou}(T_{ou}-\Omega ) & p_{ou}(t)& 
 P_{ou}\\  x_{in} 
& x_{in}(t) &\Xi_{in}= x_{in}(T_{in}+\Omega) & \Rrightarrow & 
\Xi_{ou}=x_{ou} (T_{ou}-\Omega) 
  & x_{ou}(t) & x_{ou}\\ 
 \end{array}  
\end{equation} \normalsize
where the positions are related to their \emph{celerities} 
$c_{in}(\cdot)$ and $c_{ou}(\cdot)$ by

\begin{equation} \label{e:}   
p_{in}(t) \; := \; P_{in}+ \int_{T_{in}}^{ t}c_{in}(\tau)d\tau  
\;\mbox{\rm and} \; p_{ou}(t) \; := \; P_{ou}- 
\int_{t}^{T_{ou}}c_{ou}(\tau)d\tau 
 \end{equation}
and satisfy at the junction relation

\begin{equation} \label{e:}   \begin{array}{|c|c|c|} &
\mbox{incoming dynamics} &   \mbox{outgoing dynamics}\\ \hline 
(i) &D_{in}=\varphi_{in}\Omega & D_{ou}=\varphi_{ou}\Omega 
\\   &  \Pi_{in} =
p_{in}(T_{in}+\Omega)    & \Pi_{ou} =p_{ou}(T_{ou}-\Omega )   \\
(ii) & \shortparallel   &  \shortparallel \\
 & P_{in}+ \int_{T_{in}}^{T_{in}+\Omega}c_{in}(\tau)d\tau & P_{ou} 
-\int_{T_{ou}-\Omega }^{T_{ou}}c_{ou}(\tau)d\tau\\
(iii) &\Xi_{in}=x_{in}(T_{in}+\Omega) & \Xi_{ou}=x_{ou}(T_{ou} -\Omega) \\
\\ \hline
\end{array}   \end{equation} 

\end{Definition}

\section{Transport Kernel Relations} \label{s:TransportKernelRelation}

At this stage, we assume that the junction relation $\mathcal{J}$ 
is given, and the question we investigate deals with the 
construction of transport relations $\mathcal{Q}$ containing the 
junction relation $\mathcal{J}$ and the evolution of the 
transport evolutions on the complement $\mathcal{Q} \setminus 
\mathcal{J}$ of the junction relation in the transport relation, 
which are assumed to be continuous time evolutions.  For that 
purpose, we assume that the transport evolution  is governed by a 
given   \index{transport differential inclusion} \emph{transport 
differential inclusion} controlled by incoming and outgoing 
celerities:

\begin{Definition} 
\symbol{91}\textbf{Transport Differential 
Inclusion}\symbol{93}\label{d:TransportDiffInc}\index{} Let us 
consider a set-valued map $F: (t,d,p,x) \in \mathbb{R}^{}\times 
\mathbb{R}^{}_{+} n\times K \leadsto F(t,d,p,x) \in 
\mathbb{R}^{m}$ associate with a transport state $(t,d,p,x)$ a 
subset $F(t,d,p,x)$ of monad velocities, regarded as 
 surges. We introduce also celerity functions 
$c_{in}(\cdot)$ and 
$c_{ou}(\cdot)$.\\
A \index{transport differential inclusion} \emph{transport 
differential inclusion} controls the evolution of a transport 
evolutions $t \mapsto x(t):=(x_{in} \diamondsuit_{\Omega } 
x_{ou})(t)$ of the incoming and outgoing evolutions of the monads 
\emph{outside} the junction relation regulated by celerities:

\begin{equation} \label{e:transportDiffEq} \left\{ \begin{array}{l}
\forall \; t \in [T_{in},T_{in}+\Omega ], \; \; x'_{in}(t) \; \in 
\\    \displaystyle{ F\left(t,D_{in}+ \varphi_{in}(T_{in}-t), 
P_{in}+ \int_{T_{in}}^{ t}c_{in}(\tau)d\tau,x_{in}(t) \right)}
\\  \forall \; t \in 
[T_{ou}-\Omega ,T_{ou}+\Omega ], \; \; x'_{ou}(t) \; \in\\   
\displaystyle{ F\left(t,D_{ou} + \varphi_{ou}(t-T_{ou}), P_{ou}- 
\int_{t}^{T_{ou}}c_{ou}(\tau)d\tau ,x_{ou}(t) \right)}
\end{array} \right. \end{equation}
\end{Definition}

The question is to determine when and how a  transport 
differential inclusion governs viable transport evolutions in the 
sense that

\begin{equation} \label{e:ViabTransEv} \left\{ \begin{array}{l}  
\forall \; t \in [T_{in},T_{in}+\Omega ], \; \; x_{in}(t) \; \in 
\\ \displaystyle{\mathbb{M}\left(t,D_{in}+ \varphi_{in}(T_{in}-t), 
P_{in}+ \int_{T_{in}}^{ t}c_{in}(\tau)d\tau  \right)}
\\  \forall \; t \in 
[T_{ou}-\Omega ,T_{ou}+\Omega ], \; \; x _{ou}(t) \; \in  \\
\displaystyle{ \mathbb{M}\left(t,D_{ou} + \varphi_{ou}(t-T_{ou}), 
P_{ou}- \int_{t}^{T_{ou}}c_{ou}(\tau)d\tau   \right)}
\end{array} \right. \end{equation}
which links $x_{in}$ to $x_{ou}$ from $T_{in}$ to $T_{ou}$ and 
crossing the junction on the intermodal period 
$[T_{in}+\Omega,T_{ou}-\Omega]$ during which the position 
``jumps'' from $\displaystyle{P_{in}+ 
\int_{T_{in}}^{T_{in}+\Omega}c_{in}(\tau)d\tau}$ to 
$\displaystyle{P_{ou}- \int_{T_{ou}-\Omega 
}^{T_{ou}}c_{ou}(\tau)d\tau  }$.

Hence, the question asked is: \emph{given departure states 
$(T_{in},D_{in},P_{in},x_{in})$ and arrival states 
$(T_{ou},D_{ou},P_{ou},x_{ou})$,   can they be linked by at 
least  one viable transport evolution across a junction relation 
governed by a transport differential inclusion}. 

This is a more involved version of a ``geodesic problem'' that 
can be solved thanks to the   concept of \emph{Eupalinian 
kernel}\footnote{Eupalinos, a Greek engineer, excavated around 
550 BC a 1036 m. long tunnel 180 m. below Mount Kastro for 
building an aqueduct supplying Pythagoreion (then the capital of 
Samos) with water on orders of tyrant Polycrates. \emph{He 
started to dig simultaneously the tunnel from both sides by two 
working teams who met in the center} of the channel and they had 
only 0,6 m. error. This  \emph{``Eupalinian strategy''} has been  
used ever since for building famous tunnels  or bridges by 
starting the construction at the same time from both end-points 
and proceed until they meet.}  (Section~8.5, p. 291, of 
\emph{Viability Theory.  New Directions}, \cite[Aubin, Bayen \& 
Saint-Pierre]{absp}), whereas the junction problem is a kind of 
(duration) impulse problem, dealt with in a general framework in 
Section~12.3, p.503, of this monograph, as a brief summary of an 
abundant literature (see \cite[Aubin, Lygeros, Quincampoix, 
Sastry. \& Seube]{alss99hy},\emph{Hybrid Dynamical Systems} 
\cite[Goebel et al.]{GoebelSanfeliceTeel}, Section~12.3, p.503, 
of \cite[Aubin, Bayen \& Saint-Pierre]{absp}).

\begin{Definition} 
\symbol{91}\textbf{Transport Kernel Across A 
Junction}\symbol{93}\label{}\index{} Let us consider the monad 
relation $\mathcal{M}$ describing a network, a viable junction  
relation $\mathcal{J} \subset \mathcal{M}$ of the network, 
constant incoming and outgoing fluidities $\varphi_{in}$ and 
$\varphi_{ou}$ and the differential inclusion 
(\ref{e:transportDiffEq}), p. \pageref{e:transportDiffEq} 
governing the evolution of the 
transport evolution.\\
The \index{transport kernel} \emph{transport kernel $\mbox{\rm 
Tran}_{(\ref{e:transportDiffEq})}(\mathcal{M},\mathcal{J})$ 
across a junction $\mathcal{J}$ viable in the monad relation 
$\mathbb{M}$ under the transport differential inclusion} 
(\ref{e:transportDiffEq}), p. \pageref{e:transportDiffEq} is the  
viable transport relation 
$$\mbox{\rm 
Tran}_{(\ref{e:transportDiffEq})}(\mathcal{M},\mathcal{J})  
\;\subset\; \mathcal{M}$$ made of incoming-outgoing states 
$$((T_{in},D_{in},P_{in},x_{in}), 
(T_{ou},D_{ou},P_{ou},x_{ou}))$$ which are linked by at least one 
viable \emph{transport evolution} (see 
Definition~\ref{d:TransportEvo}, p.\pageref{d:TransportEvo}) 
governed by the transport differential inclusion outside the 
junction relation $\mathcal{J}$.
  \end{Definition}

\section{Viability Construction of a Viable Transport Kernel Relation} 
\label{s:ConstructionTransKern} The definition of \emph{transport 
kernels} being set, we have to derive its properties, and, at the 
end, provide \emph{traffic regulators} indicating  mobiles 
travelling in the network  the celerities for  joining a 
departure state to an arrival states across a junction. 

For doing so, we shall construct a transport kernel as a capture 
basin of an auxiliary problem, so that the transport kernel  
inherits the properties of capture basins gathered in 
\emph{Viability Theory.  New Directions}, \cite[Aubin, Bayen \&  
Saint-Pierre]{absp}.

Let us consider  the \emph{auxiliary controlled system} defined by

\begin{equation} \label{e:JunctionCharact} \left\{ \begin{array}{llll} 
 & \mbox{incoming dynamics} & & \mbox{outgoing dynamics}\\
(i) & \overrightarrow{\tau}_{in}'(t) \; = \;+1 & (i) &  
\overleftarrow{\tau}_{ou}'(t) \; = \;  -1\\  

(ii) & \delta'_{in}(t) \; = \; - \varphi_{in}    & (ii) & 
\delta'_{ou}  \; = \; - \varphi_{ou} (t)   \\ 

(iii) & \pi'_{in}(t) \; = \;   \gamma_{in}(t) &(iii) & 
\pi'_{ou}(t) \; = \;  -\gamma_{ou}(t)
\\ \end{array} \right. 
\end{equation}
and 
\begin{equation} \label{e:}   
(iv) \; \; \xi'_{in}(t) \; \in \; 
F_{in}(\overrightarrow{\tau}_{in}(t), \delta_{in}(t), 
\pi_{in}(t), \xi_{in}(t)) \;\mbox{\rm and}\; \xi'_{ou}(t) \; \in 
\;- F_{ou}(\overleftarrow{\tau}_{ou}(t), \delta_{ou}(t), 
\pi_{ou}(t), 
\xi_{ou}(t))\\
\end{equation}
governing the evolution of pairs of   evolutions 

\begin{equation} \label{e:}   
((\overrightarrow{\tau}_{in}(t), \delta_{in}(t), \pi_{in}(t), 
\xi_{in}(t)),(\overleftarrow{\tau}_{ou}(t), \delta_{ou}(t), 
\pi_{ou}(t), \xi_{ou}(t)))
\end{equation}
  
\begin{Theorem} 
\symbol{91}\textbf{Transport Kernel 
Theorem}\symbol{93}\label{t:TransportKernelThm}\index{} Let us 
consider the monad relation $\mathcal{M}$ describing a network, a 
viable junction  relation $\mathcal{J} \subset \mathcal{M}$ of 
the network, constant incoming and outgoing fluidities 
$\varphi_{in}$ and $\varphi_{ou}$ and the differential inclusion 
(\ref{e:transportDiffEq}), p. \pageref{e:transportDiffEq} 
governing the evolution of the 
transport evolution.\\
The \index{transport kernel}  transport kernel $\mbox{\rm 
Tran}_{(\ref{e:transportDiffEq})}(\mathcal{M},\mathcal{J})$ 
across a junction $\mathcal{J}$ viable in $\mathcal{M}$ with 
under the transport differential inclusion 
(\ref{e:transportDiffEq}), p. \pageref{e:transportDiffEq} is 
equal to the capture basin $\mbox{\rm 
Capt}_{(\ref{e:JunctionCharact})}(\mathcal{M}, \mathcal{J})$ if 
the junction relation 
viable in the monad map under the auxiliary system (\ref{e:JunctionCharact}), p. \pageref{e:JunctionCharact}.\\
Therefore, the transport kernel across a junction inherits all 
properties of capture basins.
\end{Theorem}

\textbf{Proof} --- \hspace{ 2 mm} To say that 
\begin{displaymath}    
((T_{in},D_{in},P_{in},x_{in}),(T_{ou},D_{ou},P_{ou},x_{ou})) \; 
\in \; \mbox{\rm 
Capt}_{(\ref{e:JunctionCharact})}(\mathcal{M},\mbox{\rm 
Graph}(\mathcal{J}))
\end{displaymath} 
belongs to the capture basin of   junction relation $\mathcal{J}$ 
viable in  $\mathcal{M}$  means that there exit both $\Omega  
\geq 0$ (the \emph{aperture} we are looking for) and evolutions

\begin{equation} \label{e:} \left\{ \begin{array}{llll} 
 & \mbox{incoming dynamics} & & \mbox{outgoing dynamics}\\
 
(i) & t \mapsto  \overrightarrow{\tau}_{in} (t) \; = \;  T_{in} + 
t  & (i) 
 & t \mapsto  \overleftarrow{\tau}_{ou} (t) \; = \;  T_{ou} -  t 
\\  
(ii) & t \mapsto  \delta_{in}(t) \; = \; D_{in} - \varphi _{in}t 
& (ii) & t \mapsto  \delta_{ou}(t) \; = \; D_{ou} - \varphi _{ou} 
t  
\\ 
(iii) & \displaystyle{\pi_{in}(t)= P_{in}+\int_{0}^{t} 
\gamma_{in}(\tau)d\tau }
&(iii)&\displaystyle{\pi_{ou}(t) = P_{ou}- \int_{0}^{t} \gamma_{ou}(\tau)d\tau }     \\
(iv) & t \mapsto  \xi_{in}(t) &(iv) & t \mapsto  \xi_{ou}(t)   \\
\\ \end{array} \right. 
\end{equation}
starting at 
$((T_{in},D_{in},P_{in},x_{in}),(T_{ou},D_{ou},P_{ou},x_{ou}))$, 
viable in   $ \mbox{\rm Graph}({\mathbb{M}})^{2}$ such that, at 
aperture time $\Omega \geq 0$, 
\begin{equation} \label{e:} \left\{ \begin{array}{llll} 
 & \mbox{incoming dynamics} & & \mbox{outgoing dynamics}\\
 (i) & \overrightarrow{\tau}_{in} (\Omega)\; = \;  T_{in} + 
\Omega  & (i) 
 & \overleftarrow{\tau}_{ou} (\Omega)\; = \;  T_{ou} -  \Omega  
\\  
(ii) & \delta_{in}(\Omega)\; = \; D_{in} - \varphi _{in} \Omega  
\; = \; 0 & (ii) & \delta_{ou}(\Omega)\; = \; D_{ou} - \varphi 
_{ou} \Omega  \; = \; 0\\ (iii) & \displaystyle{\pi_{in}(\Omega) 
(\Omega ) = P_{in}+ \int_{0}^{\Omega 
}\gamma_{in}(\tau)d\tau}&(iii) &  \displaystyle{\pi_{ou}(\Omega) 
(\Omega )P_{ou}-\int_{0}^{\Omega }\gamma_{ou}(\tau)d\tau } \\
(iv) & \xi_{in}(\Omega)&(iv) & \xi_{ou}(\Omega)   \\
 \\ \end{array} 
\right. 
\end{equation}
the pair 
$$(( T_{in} + \Omega,D_{in}-\varphi _{in}\Omega 
,\pi_{in}(\Omega),\xi_{in}(\Omega)),
 (T_{ou}-\Omega,D_{ou}-\varphi_{ou} \Omega, \pi_{ou}(\Omega), 
\xi_{ou}(\Omega)))$$ belongs to the junction relation $ 
\mathcal{J}$.

By definition of the junction  relation $ \mathcal{J} $, this 
implies that   $D_{in}= \varphi_{in} \Omega $ and $D_{ou}= 
\varphi _{ou} \Omega $. Consequently, the ratios 
\begin{equation} \label{e:}    
\frac{D_{in}}{\varphi_{in}}  \; := \;  
\frac{D_{ou}}{\varphi_{ou}} \; = \; \Omega 
\end{equation} 
are all equal to the aperture $\Omega $.

By making the changes of variables $t \mapsto t-T_{in}$ in the 
evolutions of the auxiliary incoming states and $t \mapsto T_{ou} 
-t$ in the auxiliary outgoing states, we observe that 

\begin{enumerate}  
\item  $\overrightarrow{\tau}_{in}(t-T_{in})=t$ and 
$\tau_{ou}(T_{ou}-t) =   t$;

\item $d_{in}(t):=\delta_{in}(t-T_{in})=D_{in} + 
\varphi_{in}T_{in} -\varphi_{i} t$ and 
$d_{ou}(t)=\delta_{ou}(T_{ou} -t)=D_{ou} - \varphi_{ou}T_{ou} 
+\varphi_{ou} t$;

\item defining celerities $c_{in}(t):=\gamma_{in}(t-T_{in})$ and  
$c_{ou}(t):=\gamma_{ou}(T_{ou}-t)$,  

\begin{equation} \label{e:} \left\{ \begin{array}{l}  
\displaystyle{p_{in}(t) \; := \; \pi_{in}(t-T_{in}) \; = \; 
P_{in}+\int_{T_{in}}^{t}c_{in}(\tau)d\tau } \\ 
\displaystyle{p_{ou}(t)\; := \; \pi_{ou}(T_{ou} -t) \; = \; 
P_{ou}- \int_{t}^{T_{ou}}c_{ou}(\tau)d\tau}
\end{array} \right. 
\end{equation}

\item $x_{in}(t):=\xi_{in}(t-T_{in}) $ and 
$x_{ou}(t):=\xi_{ou}(T_{ou} -t) $, satisfying the extremal 
conditions $\xi_{in}(T_{in}):=\xi_{in} $ and $ 
\xi_{ou}(T_{ou}):=\xi_{ou} $,

\end{enumerate}
and the requirements for describing a \index{transport evolution} 
\emph{transport evolution} starting at 
$(T_{in},D_{in},P_{in},x_{in})$ and arriving at 
$(T_{ou},D_{ou},P_{ou},x_{ou})$ across the junction $\mathcal{J}$ 
with aperture $\Omega \geq 0$ according to 
Definition~\ref{d:TransportEvo}, p.\pageref{d:TransportEvo}. 
\hfill $\;\; \blacksquare$ \vspace{ 5 mm}
 
 \textbf{Remark} --- \hspace{ 2 mm}  If the junction relation $ 
 \mathcal{J}  := 
\{(\Sigma_{in},\Sigma_{ou},\Pi_{in},\Pi_{ou})\}$ is reduced to a 
singleton, this implies $T_{in} \; = \; \Sigma_{in}-\Omega $ and 
$T_{ou} \; = \; \Sigma_{ou}+\Omega $ as well as 
$\displaystyle{P_{in}=\Pi_{in}-\int_{0}^{\Omega} 
c_{in}(\tau)d\tau}$ and $\displaystyle{P_{ou}=\Pi_{ou} 
-\int_{T_{ou}-\Omega }^{T_{ou}} c_{ou}(\tau)d\tau}$. In 
particular, if the junction is an time impulse junction at 
impulse date $\Sigma$, then we infer that the initial and arrival 
times are equal to  $T_{in} \; = \; \Sigma -\Omega $ and $T_{ou} 
\; = \; \Sigma +\Omega $. \hfill $\;\; \blacksquare$ \vspace{ 5 
mm}

We  introduce supplementary ``decomposition'' condition: 

\begin{Definition} 
\symbol{91}\textbf{Decomposable Transport Relations Outside a 
Junction Relation}\symbol{93}\label{}\index{}  A transport 
relation is \index{decomposable transport relations outside a 
junction relation} \emph{decomposable outside} a junction 
relation if complement  
\begin{equation} \label{e:} \mathcal{Q} \setminus 
\mathcal{J}\; = \; \mathcal{Q}_{in} \times \mathcal{Q}_{ou}  
\end{equation}  of the junction relation $\mathcal{J}$ in the viable transport 
relation $\mathcal{Q}$ is split as the product of  disjoint 
prejunction transport relation $\mathcal{Q}_{in} =: \mbox{\rm 
Graph}(\mathbb{M}_{\mathcal{Q}_{in}})$ and of a postjunction 
transport relation $\mathcal{Q}_{ou}=: \mbox{\rm 
Graph}(\mathbb{M}_{\mathcal{Q}_{ou}})$, where
\begin{equation} \label{e:} \left\{ \begin{array}{l}    
\mathcal{Q}_{in} := \{(t ,d ,p  ,x ) \in \mathcal{Q} \; \mbox{ 
such that} \; t \; < \;  \Sigma^{\mathcal{J}}_{in}  := 
\inf_{((\Sigma_{in},\Pi_{in},\Xi_{in}),(\Sigma_{ou}, 
\Pi_{ou},\Xi_{ou}))}\Sigma_{in}\}
\\
\mathcal{Q}_{ou} := \{(t ,d ,p  ,x ) \in \mathcal{Q} \; \mbox{ 
such that} \; t \;> \; \Sigma^{\mathcal{J}}_{ou}  := 
\inf_{((\Sigma_{in},\Pi_{in},\Xi_{in}),(\Sigma_{ou}, 
\Pi_{ou},\Xi_{ou}))}\Sigma_{ou}\}
\end{array} \right. 
\end{equation}
so that  $((T_{in},D_{in},P_{in},x_{in}), 
(T_{ou},D_{ou},P_{ou},x_{ou})) \in \mathcal{Q} \setminus 
\mathcal{J}$ if and only if 

\begin{equation} \label{e:} \left\{ \begin{array}{l}  
\forall \; t \in  \symbol{91}T_{in}, T_{in}+\Omega \symbol{91}, 
\; \;x_{in}(t) \in \mathbb{M}^{\mathcal{J}}_{in}(t, 
D_{in}+\varphi_{in}(t-T_{in}), p_{in}(t)) \\
\forall \; t \in  \symbol{93}T_{ou}- \Omega , T_{ou}\symbol{93}, 
\; \;x_{ou}(t) \in \mathbb{M}^{\mathcal{J}}_{ou}(t, 
D_{ou}-\varphi_{ou}( T_{ou}-t), p_{ou}(t))
\end{array} \right.
\end{equation}

 \end{Definition}

Examples of decomposable transport relations are provided by 
monad relations satisfying a ``safety condition'': 

\begin{Definition} 
\symbol{91}\textbf{Safety 
Condition}\symbol{93}\label{d:safety}\index{temporal profile} We 
associate with the monad  relation $\mathcal{M}$ its 
\index{safety condition} \emph{temporal profile} $P_{\mathbb{M}}: 
\mathbb{R}^{}  \leadsto \mathbb{R}^{}_{+} \times K \times 
\mathbb{R}^{m}$ defined by 
\begin{equation} \label{e:}   
(d,p,x) \; \in \; P_{\mathbb{M}}(t) \;\mbox{\rm if and only if}\; 
x \; \in \; \mathbb{M}(t,d,p)
\end{equation}
The monad  relation $\mathcal{M}$ satisfies the \emph{safety 
condition} if whenever $t \ne s$, then $P_{\mathbb{M}}(t) \cap 
P_{\mathbb{M}}(s) = \emptyset$.
\end{Definition}

We observe that

\begin{Lemma} 
\symbol{91}\textbf{Safe Monad Relations are 
Decomposable}\symbol{93}\label{}\index{} Any  monad  relation 
$\mathcal{M}$ satisfying the  safety condition is decomposable 
outside any junction relation $\mathcal{J} \subset \mathcal{M}$.
\end{Lemma}
 
\textbf{Proof} --- \hspace{ 2 mm} Indeed, let 
$((t_{in},d_{in},p_{in},x_{in}),(t_{ou},d_{ou},p_{ou},x_{ou}))$ 
belong to the transport relation $\mathcal{Q}$. Since it does not 
belong to the junction relation $\mathcal{J}$, then $t_{in} \leq  
 \Sigma^{\mathcal{J}}_{in}  \leq 
\Sigma^{\mathcal{J}}_{ou} \leq   t_{ou}$ and $ t_{in} \ne  
t_{ou}$, since either $d_{in}>0$ and thus, $t_{in} < 
\Sigma^{\mathcal{J}}_{in}$ or $d_{ou}>0$ and thus, $t_{ou} > 
\Sigma^{\mathcal{J}}_{ou}$. Consequently,   \emph{temporal 
profiles} $P_{\mathbb{M}}(t_{in}) \cap P_{\mathbb{M}}(t_{ou}) = 
\emptyset$ as subsets of $\mathbb{R}^{}_{+} \times K \times 
\mathbb{R}^{m}$.   \hfill $\;\; \blacksquare$ \vspace{ 5 mm}

\section{Construction of Transport Regulators} \label{s:TransportRegulator}
 
The question arises whether one can construct   viable transport 
relation $\mathcal{Q}$. Viability theorems allow us to derive 
some properties of transport kernel relation
 $\mbox{\rm 
Tran}_{(\ref{e:transportDiffEq})}(\mathcal{M},\mathcal{J})$ 
across a junction $\mathcal{J}$ viable in the monad relation 
$\mathbb{M}$ under the transport differential inclusion 
(\ref{e:transportDiffEq}), p. \pageref{e:transportDiffEq}.   

\begin{Theorem} 
\symbol{91}\textbf{Properties of the Transport Kernel 
Relation}\symbol{93}\label{}\index{} Assume that the graphs of 
the monad relation $\mathcal{M}$  and of the junction relation 
$\mathcal{J}$ are closed,  that  incoming and outgoing fluidities 
$\varphi_{in}$ and $\varphi_{ou}$ are constant  and that the 
set-valued map  $F: (t,d,p,x) \in \mathbb{R}^{}\times 
\mathbb{R}^{}_{+}  \times K \leadsto F(t,d,p,x) \in 
\mathbb{R}^{m}$ (see Definition~\ref{d:TransportDiffInc}, 
p.\pageref{d:TransportDiffInc}) is Marchaud. \\Then the transport 
kernel relation $\mbox{\rm 
Tran}_{(\ref{e:transportDiffEq})}(\mathcal{M},\mathcal{J})$ 
across a junction $\mathcal{J}$ viable in the monad relation 
$\mathbb{M}$ under the transport differential inclusion  
(\ref{e:transportDiffEq}), p. \pageref{e:transportDiffEq},
\begin{enumerate}   
\item is closed,

\item if the transport kernel is not empty and decomposable 
outside the junction, for all departure and arrival states 
$(T_{in},D_{in},P_{in},x_{in})$ and $ 
(T_{ou},D_{ou},P_{ou},x_{ou}) \in \mbox{\rm 
Tran}_{(\ref{e:transportDiffEq})}(\mathcal{M},\mathcal{J}) $, 
then  all transport evolutions governed by 
(\ref{e:transportDiffEq}) linking them and which are  viable in 
monad relation $\mathcal{M}$ in the sense of 
(\ref{e:ViabTransEv}), p. \pageref{e:ViabTransEv}: 
\begin{displaymath}    
\left\{ \begin{array}{l}  \forall \; t \in [T_{in},T_{in}+\Omega 
], \; \; x_{in}(t) \; \in 
\\ \displaystyle{\mathbb{M}\left(t,D_{in}+ \varphi_{in}(T_{in}-t), 
P_{in}+ \int_{T_{in}}^{ t}c_{in}(\tau)d\tau  \right)}
\\  \forall \; t \in 
[T_{ou}-\Omega ,T_{ou}+\Omega ], \; \; x _{ou}(t) \; \in  \\
\displaystyle{ \mathbb{M}\left(t,D_{ou} + \varphi_{ou}(t-T_{ou}), 
P_{ou}- \int_{t}^{T_{ou}}c_{ou}(\tau)d\tau   \right)}
\end{array} \right.
\end{displaymath}
are actually viable in the complement $\mbox{\rm 
Tran}_{(\ref{e:transportDiffEq})}(\mathcal{M},\mathcal{J}) 
\setminus \mathcal{J}$ of the junction map in the sense that

\begin{equation} \label{e:TransKernViabPr}    
\left\{ \begin{array}{l}  \forall \; t \in [T_{in},T_{in}+\Omega 
], \; \; x_{in}(t) \; \in 
\\ \displaystyle{\mathbb{M}^{\mathcal{J}}_{in}\left(t,D_{in}+ \varphi_{in}(T_{in}-t), 
P_{in}+ \int_{T_{in}}^{ t}c_{in}(\tau)d\tau  \right)}
\\  \forall \; t \in 
[T_{ou}-\Omega ,T_{ou}+\Omega ], \; \; x _{ou}(t) \; \in  \\
\displaystyle{\mathbb{M}^{\mathcal{J}}_{ou}\left(t,D_{ou} + 
\varphi_{ou}(t-T_{ou}), P_{ou}- 
\int_{t}^{T_{ou}}c_{ou}(\tau)d\tau   \right)}
\end{array} \right.
\end{equation} 
\end{enumerate}\end{Theorem}

The necessary condition of the Viability Theorem implies that 
whenever $F$ is Marchaud, the viable transport evolutions $t 
\mapsto ((t,d_{in}(t),p_{in}(t),x_{in}(t)), 
(t,d_{ou}(t),p_{ou}(t),x_{ou}(t)))$ governed by differential 
inclusions (\ref{e:transportDiffEq}), p. 
\pageref{e:transportDiffEq}, and viable in the sense of 
(\ref{e:TransKernViabPr}), p. \pageref{e:TransKernViabPr}: For 
that purpose, we have to introduce the (convexified) tangent 
cones to the partial transport relations.

\begin{Definition} 
\symbol{91}\textbf{Transport Regulators}\symbol{93} 
\label{d:TransRegul} \index{transport regulator} The transport 
regulators $R^{(\mathcal{M},\mathcal{J})}_{in}$ and 
$R^{(\mathcal{M},\mathcal{J})}_{ou}$  are defined by

\begin{equation} \label{e:FormTransReg} \left\{ \begin{array}{l}  
c_{in} \; \in \;  R^{(\mathcal{M},\mathcal{J})}_{in}(t,d,p,x) \;\mbox{\rm if and only if }\\
(1,-\varphi_{in},c_{in},F(t,d,p,x)) \cap 
T_{Q^{(\mathcal{M},\mathcal{J})}_{in}}(t,d,p,x) \ne 
\emptyset \\
c_{ou} \; \in \;  R^{(\mathcal{M},\mathcal{J})}_{in} (t,d,p,x) \;\mbox{\rm if and only if }\\
(1,+\varphi_{ou},c_{ou},F(t,d,p,x)) \cap 
T_{Q^{(\mathcal{M},\mathcal{J})}_{ou}}(t,d,p,x) \ne 
\emptyset \\
\end{array} \right. 
\end{equation}

\end{Definition}

We thus infer that the velocities sent to the mobiles of the 
network crossing one junction are controlled by the transport 
regulators. 
 
\begin{Theorem} 
\symbol{91}\textbf{Transport 
Regulations}\symbol{93}\label{}\index{}  Assume that the monad 
relation $\mathcal{M}$  and of   junction relation $\mathcal{J}$ 
are closed, that $\mathcal{M}$ is decomposable outside 
$\mathcal{J}$, that  incoming and outgoing fluidities 
$\varphi_{in}$ and $\varphi_{ou}$ are constant  and that the 
set-valued map  $F: (t,d,p,x) \in \mathbb{R}^{}\times 
\mathbb{R}^{}_{+}  \times K \leadsto F(t,d,p,x) \in 
\mathbb{R}^{m}$ (see Definition~\ref{d:TransportDiffInc}, 
p.\pageref{d:TransportDiffInc}) is Marchaud. Then the   transport 
evolutions in the transport kernel relation are controlled by the 
celerities provided by the transport regulators:

 \begin{equation} \label{e:} \left\{ \begin{array}{l}  
c_{in} (t)\; \in \;  R^{(\mathcal{M},\mathcal{J})}_{in} 
(t,D_{in}+\varphi_{in}(t-T_{in}),
 p_{in}(t),x_{in}(t))   \\
c_{ou}(t)
 \; \in \;  R^{(\mathcal{M},\mathcal{J})}_{ou} (t, 
D_{ou}-\varphi_{ou}( T_{ou}-t), p_{ou}(t),x_{ou}(t)) \\
\end{array} \right. 
\end{equation}
\end{Theorem}

\textbf{Remark} --- \hspace{ 2 mm} If we regard the incoming and 
outgoing relations $Q^{(\mathcal{M},\mathcal{J})}_{in}$ and 
$Q^{(\mathcal{M},\mathcal{J})}_{ou}$ as the graphs of the monad 
maps $ \mathcal{M}_{in}$ and $ \mathcal{M}_{ou}$, and recalling 
that the tangent cone to the graph of a set-valued map is the 
graph of its derivative: 
 
\begin{equation} \label{e:} \left\{ \begin{array}{l}  
T_{Q^{(\mathcal{M},\mathcal{J})}_{in}}(t,d,p,x) \; := \; 
\mbox{\rm Graph}(D\mathbb{M}^{\mathcal{J}}_{in}(t,d,p,x))\\
T_{Q^{(\mathcal{M},\mathcal{J})}_{ou}}(t,d,p,x) \; := \; 
\mbox{\rm Graph}(D\mathbb{M}^{\mathcal{J}}_{ou}(t,d,p,x))
\end{array} \right. \end{equation}
formulas (\ref{e:FormTransReg}), p. \pageref{e:FormTransReg} 
defining the transport regulators   
$R^{(\mathcal{M},\mathcal{J})}_{in}$ and 
$R^{(\mathcal{M},\mathcal{J})}_{ou}$ can be rewritten in 
``differential form''

\begin{equation} \label{e:FormTransRegdif} \left\{ \begin{array}{l}  
c_{in} \; \in \;  R^{(\mathcal{M},\mathcal{J})}_{in}(t,d,p,x) 
\;\mbox{\rm if and only if}  \\  F(t,d,p,x)\cap 
D\mathbb{M}^{\mathcal{J}}_{in}(t,d,p,x)
(1,-\varphi_{in}, c_{in}) \; \ne \; \emptyset\\
c_{ou} \; \in \;  R^{(\mathcal{M},\mathcal{J})}_{ou}(t,d,p,x) 
\;\mbox{\rm if and only if} \\  F(t,d,p,x)\cap 
D\mathbb{M}^{\mathcal{J}}_{ou}(t,d,p,x) (1,\varphi_{ou}, c_{ou}) 
\; \ne \; \emptyset\\ 
\end{array} \right. 
\end{equation}
Except for particular cases, the formulations in terms of 
derivatives of set-valued maps to obtain set-valued versions of 
conservation laws or non smooth solutions of 
Hamilton-Jacobi-Bellman is not always useful. \hfill $\;\; 
\blacksquare$ \vspace{ 5 mm}

\clearpage

\tableofcontents


\begin{thebibliography}{a}
 \bibitem{Anita} Anita S., (2000)
\textbf{\emph{Analysis and Control of Age-Dependent population 
Dynamics}}, Kluwer Academic Publishers
  
 


\bibitem{avt} Aubin J.-P.    (1991)
\textbf{\emph{Viability Theory}}, Birkh\"auser

\bibitem{a92ia} Aubin J.-P.    (1996)
\textbf{Neural Networks and Qualitative Physics: a Viability 
Approach}, Cambridge University Press

\bibitem{aconcom98} Aubin J.-P.    (1998)
\emph{Connectionist Complexity and its Evolution}, in \textbf{ 
Equations aux d\'eriv\'ees partielles, Articles d\'edi\'es \`a 
J.-L. Lions}, 50-79, Elsevier
 

\bibitem{abb99evs} Aubin J.-P., Bonneuil N. and Maurin F. (2000)
Non-linear Structured Population Dynamics with Co-Variables, 
\emph{Mathematical Population Studies}, 91, 1-31

\bibitem{a01ctfcr} Aubin J.-P.    (2003)
\emph{Regulation of the Evolution of the Architecture of a 
Network by Connectionist Tensors Operating on Coalitions of 
Actors},J. Evolutionary Economics, 13,95-124

\bibitem{mded} Aubin J.-P.    (2010)
\textbf{La  mort  du  devin,  l'émergence du  démiurge. Essai sur 
la contingence, la viabilité et l'inertie des systèmes}, Éditions 
Beauchesne

\bibitem{AUB-Leit09}  Aubin J.-P. (2010)
\emph{Macroscopic Traffic Models: Shifting from Densities to 
``Celerities''}, Applied Mathematics and Computation, 217, 
963-971, \url{http://dx.doi.org/10.1016/j.amc.2010.02.032}

\bibitem{RegulationBirths-11}  Aubin J.-P.
(2011), \emph{Regulation of Births for Viability of Populations 
Governed by  Age-Structured Problems}, Journal of Evolution 
Equations, (DOI: 
\url{http://dx.doi.org/10.1007/s00028-011-0125-z})

\bibitem{AUB-10-STR}  Aubin J.-P. (2011)
\emph{Viability Solutions to Structured Hamilton-Jacobi Equations 
under Constraints}, SIAM Journal of Control and Optimization, 49, 
1881-1915,  \url{http://dx.doi.org/10.1137/10079567X}

\bibitem{Lax-Hopf}  Aubin J.-P.
(2012), \emph{Lax-Hopf Formula and Max-Plus Properties of 
Solutions to Hamilton-Jacobi equations}, NoDEA, DOI: (DOI) 
10.1007/s00030-012-0188-8

\bibitem{aub-12-Durance} Aubin J.-P. (2013)
\emph{Chaperoning State Evolutions by Variable Durations}, SIAM 
Journal of Control and Optimization, DOI. 10.1137/120879853 

\bibitem{TM} Aubin J.-P. (2013) 
\emph{Time and Money. How Long  and How Much Money is Needed to 
Regulate a Viable  Economy}, Springer

\bibitem{Trend} Aubin J.-P. (submitted) 
\emph{Retro-Prospective Differential Inclusions and their Control 
by the Differential Connection Tensors of their Evolutions: The 
trendometer},



\bibitem{absp} Aubin J.-P., Bayen A. \& Saint-Pierre P. (2011)
\textbf{\emph{Viability Theory.  New Directions}}, Springer



\bibitem{absp08hj} Aubin J.-P., Bayen A. \& Saint-Pierre P. (2008)
\emph{Dirichlet Problems for some Hamilton-Jacobi Equations With 
Inequality Constraints}, SIAM Journal of Control and 
Optimization, 47, 2348-2380, DOI: 
\url{http://dx.doi.org/10.1137/060659569}


\bibitem{JPA-LC-Cournot-13} Aubin J.-P. and Chen Luxi. (in preparation) 
Cournot Maps for Intercepting Evader Evolutions by a Pursuer 

\bibitem{JPA-LC-Lax-13} Aubin J.-P. and Chen Luxi. (in preparation) 
\emph{Generalized Lax-Hopf formulas for Cournot Maps and 
Hamilton-Jacobi-McKendrik equations},

\bibitem{ADA-Traffic} Aubin J.-P. \& Désilles A. (in preparation) 
\textbf{\emph{Mathematical Approaches to Traffic Management}}, 
Springer


\bibitem{af90sva} Aubin J.-P. et Frankowska H. (1990) 
\textbf{\emph{Set-valued analysis}}, Birkhäuser

\bibitem{alss99hy} Aubin J.-P., Lygeros J., Quincampoix M., Sastry S. \& Seube
N.    (2002) \emph{Impulse Differential Inclusions: A Viability 
Approach to Hybrid Systems}, IEEE Transactions on Automatic 
Control, 47, 2-20
 
\bibitem{Aub-Mart08} Aubin J.-P. \& S.  Martin S. (2009)
Travel Time Tubes Regulating Transportation Traffic, 
\emph{Contemporary mathematics},   1-25


\bibitem{BCSP07a}
Bayen A., Claudel C.  Saint-Pierre P. (2007) Viability-based 
computations of solutions to the   {H}amilton-{J}acobi-{B}ellman 
equation,  In  {\em Hybrid Systems: Computation and Control}, 
Lecture Notes in Computer Science 4416, pp 645-649. 
Springer-Verlag 



\bibitem{BCSP07b}
Bayen A., Claudel C.  Saint-Pierre P. (2007) Computations of 
solutions to the Moskowitz Hamilton-Jacobi-Bellman equation under 
viability constraints, Proceedings of the \textit{46th IEEE 
Conference on Decision and Control (CDC)}, New Orleans, LA, 
available online at 
\url{http://www.ce.berkeley.edu/$\sim$bayen/conferences/cdc07a.pdf} 

\bibitem{Chagnon} Chagnon N. (2013) \emph{Noble Savages}, Simon \& Schuster

 
\bibitem{CB08a}
  Claudel C. and  Bayen A. (2008)  Solutions to switched 
{H}amilton-{J}acobi equations and conservation laws using hybrid 
components,  in  {\em Hybrid Systems: Computation and Control}, 
Lecture Notes in Computer Science (Egerstedt M and Mishra B., 
editor), 4981, pp 101-115. Springer-Verlag,  


\bibitem{CB08b}
 Claudel C. and  Bayen A. (2008) Computations of solutions to the 
Moskowitz Hamilton-Jacobi-Bellman equation under trajectory 
constraints,\textit{46th IEEE Conference on Decision and Control 
(CDC)}, Cancun, Mexico.

 
\bibitem{DaganzoCell94}
Daganzo C. (1994)
 The cell transmission model: a dynamic representation of highway
  traffic consistent with the hydrodynamic theory.
 {\em Transportation Research}, 28B(4):269--287 
\bibitem{DaganzoCell95}
Daganzo C. (1995)
 The cell transmission model, part {II}: network traffic.
 {\em Transportation Research}, 29B(2):79--93   

\bibitem{DaganzoCell95}
Daganzo C. (1995)
 The cell transmission model, part {II}: network traffic.
 {\em Transportation Research}, 29B(2):79--93   

\bibitem{DaganzoVP-1}
Daganzo C.
 A variational Formulation of Kinematic Wages: Solution Method
 {\em Transportation Research}, 39B: 934--950, (2005)

\bibitem{DaganzoVP-2}
Daganzo C. (2005)
 A variational Formulation of Kinematic Wages: Basic Theory and Complex Boundary Conditions
 {\em Transportation Research}, 39B: 187--196  

\bibitem{Anya} Désilles A. (2013) Viability approach to 
Hamilton-Jacobi-Moskowitz problem involving variable regulation 
parameters,  Applied Math Journal Networks and Heterogeneous Media
 
   
\bibitem{Evans98}
 Evans L. C. (1998)
 {\em Partial Differential Equations}.
 American Mathematical Society, Providence, RI, 

 
\bibitem{fh91cdc}
 Frankowska H. (1991) \emph{Lower semicontinuous solutions to 
Hamilton-Jacobi-Bellman equations}, in 
\textbf{\textit{Proceedings of the 30th IEEE Conference on 
Decision and Control}}, Brighton, UK, 


\bibitem{HJB92}
 Frankowska H. (1993) \emph{Lower semicontinuous solutions of 
Hamilton-Jacobi-Bellman equations}, SIAM J. Control Optim., 31 
  257-272.
  
\bibitem{GoebelSanfeliceTeel} Goebel R., Sanfelice R. G \& Teel A. R.
(2012) \emph{\emph{Hybrid Dynamical Systems. Modeling, Stability 
and Robustness}}, Princeton University Press
 
\bibitem{Gould77} Gould S. J. (1977) \emph{Phylogenesis and Ontogenesis}, Harvard University Press
 

\bibitem{Greenshields} 
Greenshields B. D.  (1933) The Photographic Method of studying 
Traffic Behaviour,  in \emph{Proc. of the 13 th Annual Meeting of 
the Highway Research Board} 

\bibitem{im95aspd} Iannelli M. (1995)
\textbf{\emph{Mathematical theory of travel-time-structured 
traffic management}}, Giardini, Pisa

\bibitem{ke297dem} Keyfitz B. \& Keyfitz N.   (1997)
        {\it The McKendrick partial differential equation and its uses in
epidmiology and population study}, Math. Comput. Modelling, 26, 
1-9

\bibitem{mck26dem} Mckendrick A.G. (1926)
Applications of mathematics to medical problems, \emph{Proc. 
Edinb. Math. Soc}, 44, 98-130



\bibitem{RockNetw} Rockafellar R.T. (1984) \emph{Network flows and monotropic 
optimization}, Wiley

\bibitem{rw91nsa} Rockafellar R.T. \& Wets R. (1997) 
\textbf{\emph{Variational Analysis}}, Springer-Verlag

\bibitem{Villani} Villani C. (2008) \emph{Optimal Transport: Old and 
New}, Springer
 
\bibitem{Villani2} Villani C. (2010) \emph{Optimal Transport}, Springer
 
 
\bibitem{foe59dem} Von Foerster H. (1959)
Some remarks on changing populations, In \emph{The kinetics of 
cell proliferation}, 382-407

\bibitem{web85age} Webb G. (1985)
\textbf{\emph{Theory of nonlinear age-dependent population 
dynamics}}, Marcell Dekker
 
\end{thebibliography}
\end{document}